\documentclass[12pt]{amsart}

\usepackage{amsrefs}
\usepackage{amssymb}
\usepackage{hyperref}
\usepackage{fancyhdr}
\usepackage{xcolor}
\usepackage{bbm}

\newtheorem{theorem}{Theorem}
\newtheorem*{theorem*}{Theorem}

\newtheorem{claim}[theorem]{Claim}

\newtheorem*{remark*}{Remark}

\newtheorem{maintheorem}{Theorem}

\theoremstyle{definition}

\newtheorem*{definition*}{Definition}
\newtheorem*{lemma*}{Lemma}
\usepackage{graphicx}
\usepackage{pstricks, enumerate, pst-node, pst-text, pst-plot}

\numberwithin{equation}{section}
\numberwithin{theorem}{section}

\newcommand{\N}{\mathbb{N}}

\newcommand{\Z}{\mathbb{Z}}

\def\cG{{\mathcal G}}
\def\bT{{\mathbb T}}

\begin{document}

\title[]{Transitive graphs uniquely determined by their local structure}

\author{Joshua Frisch and Omer Tamuz}
\address{Department of Mathematics, Massachusetts Institute of
  Technology, Cambridge MA 02139, USA.}


\thanks{J.\ Frisch was supported by MIT's Undergraduate Research
  Opportunities Program. This research was partially conducted at
  Microsoft Research, New England.}

\date{\today}

\begin{abstract}
  We show that the ``grandfather graph'' has the following property:
  it is the unique completion to a transitive graph of a large enough
  finite subgraph of itself.
\end{abstract}

\maketitle

\section{Introduction}
Transitive graphs ``look the same from the point of view of every
vertex''; all vertices play the same role in their geometry. Thus they
are a natural model for a discrete, homogeneous geometrical space. In
this paper we study transitive graphs whose local structure determines
their global structure.

Formally, let $\cG$ be the set of finite or countably infinite,
simple, undirected, locally finite, connected, vertex transitive
graphs. Given a $G \in \cG$ and an $r \in \N$, a ball of radius $r$ in
$G$ is the subgraph induced by all vertices at distance at most $r$
from some vertex in $G$.  We say that $G=(V,E) \in \cG$ is {\em
  isolated} if it has the following property: there exists an $r \in
\N$ large enough so that, if a ball of radius $r$ in some $H \in \cG$
is isomorphic to the ball of radius $r$ in $G$, then $H$ is isomorphic
to $G$. Intuitively, the structure of the ball of radius $r$ in $G$
determines $G$ uniquely.

Clearly, every finite transitive graph is isolated: one can take $r$
to be the radius of $G$. However, it is not obvious that there are any
{\em infinite} graphs that have this property. For example, let $G$ be
the bi-infinite chain (i.e., the Cayley graph of $\Z$). Then any ball
in $G$ can be completed to a large enough finite chain. In this paper
we give an example of an isolated infinite graph, namely Trofimov's
{\em grandfather graph}~\cite{trofimov1985automorphism}. In fact, we
give a countable family of such examples.

Note that the grandfather graph is not unimodular, and so cannot
locally resemble finite graphs. The novelty is therefore that it can
also not locally resemble any other infinite graph. It would be
interesting to find an example of a isolated, finite, unimodular
graph.

This question can be formulated as one of finding isolated points in a
natural topology on the set of transitive graphs, namely the
Benjamini-Schramm topology~\cite{benjamini2001recurrence,
  aldous2003objective}. A number of interesting questions arise: what
is the Cantor-Bendixson rank of this space? Which graphs are left
after the isolated points are repeatedly removed? And what generic
properties does such a graph have?

These and similar questions have been previously addressed in regard
to the related space of {\em marked
  groups}~\cite{champetier2000espace, cornulier2009cantor}. In
particular, Cornulier, Guyot and Pitsch~\cite{de2007isolated}
characterize the isolated points in that space. It would be
interesting to understand if the (unlabeled) Cayley graphs of these
groups are isolated in the space of transitive graphs.

\subsection{Acknowledgments}
We would like to thank Russell Lyons for helpful discussions.

\section{Formal definitions and results}

\subsection{Transitive graphs}
Let $G =(V,E)$ be a graph. We will study the set of graphs  with
the following properties:
\begin{itemize}
\item $V$ is finite or countably infinite.
\item $G$ is simple and undirected: $E$ is a symmetric relation on
  $V$.
\item $G$ is locally finite: the number of edges incident on each
  vertex is finite.
\item $G$ is connected: there is a path between every pair of
  vertices.
\item $G$ is vertex transitive; we next define this notion.
\end{itemize}

A {\em graph isomorphism} between $G=(V,E)$ and $H=(U,F)$ is a
bijection $h \colon V \to U$ such that $(u,w) \in E$ if and only if
$(h(u),h(w)) \in F$.  A graph automorphism is a graph isomorphism from
a graph to itself. A graph $G=(V,E)$ is said to be {\em vertex
  transitive} if its automorphism group acts transitively on its
vertices. That is, if for every $u,w \in V$ there exists an
automorphism $h$ such that $h(u)=w$. The {\em isomorphism class} of a
transitive graph $G$ is the set of graphs $H$ that are isomorphic to
$G$.  We denote by $\cG$ the set of isomorphism classes of graphs with
the properties described above. In this paper, we will, whenever
unambiguous, refer to ``graph isomorphism classes'' simply as
``graphs'', and likewise simply denote by $G$ the isomorphism class of
$G$. We will accordingly write $G=H$ whenever $G$ and $H$ are in the
same isomorphism class.

Given $G=(V,E)\in \cG$ and $r \in \N$, let $B_r(G) = (V_r,E_r)$ be the
ball of radius $r$ in $G$. This is the finite induced subgraph of $G$
whose vertices $V_r$ are all the vertices at distance at most $r$ from
some vertex of $G$, and whose edges $E_r$ are the edges of $G$ whose
vertices are both in $V_r$. Since we are concerned with graph
isomorphism classes, and since $G$ is vertex transitive, it does not
matter with which vertex of $G$ we choose to construct $B_r(G)$.

\subsection{The Benjamin-Schramm topology and isolated points}
The Benjamini-Schramm topology~\cite{benjamini2001recurrence,
  aldous2003objective} on $\cG$ is defined by the following
metric. Given $G,H \in \cG$, let
\begin{align*}
  D(G,H) = \sup\{2^{-r}\,:\,B_r(G) = B_r(H)\}.
\end{align*}
It is straightforward to verify that this is indeed a metric. In fact,
this topology is Polish and zero-dimensional. The sets $\cG_d$
consisting of the graphs with degree $d$ are compact in this topology.

We say that $G \in \cG$ is {\em isolated} if it is an isolated point
in this topology. By the above definition, this means that there
exists an $r \in \N$ such that whenever $B_r(G)=B_r(H)$ then
$G=H$. Since $B_r(G)=G$ for every finite $G$ and $r$ large enough, it
follows immediately that all the finite graphs are isolated.

\subsection{The grandfather graph}

\begin{figure}[h]
  \centering
  \includegraphics{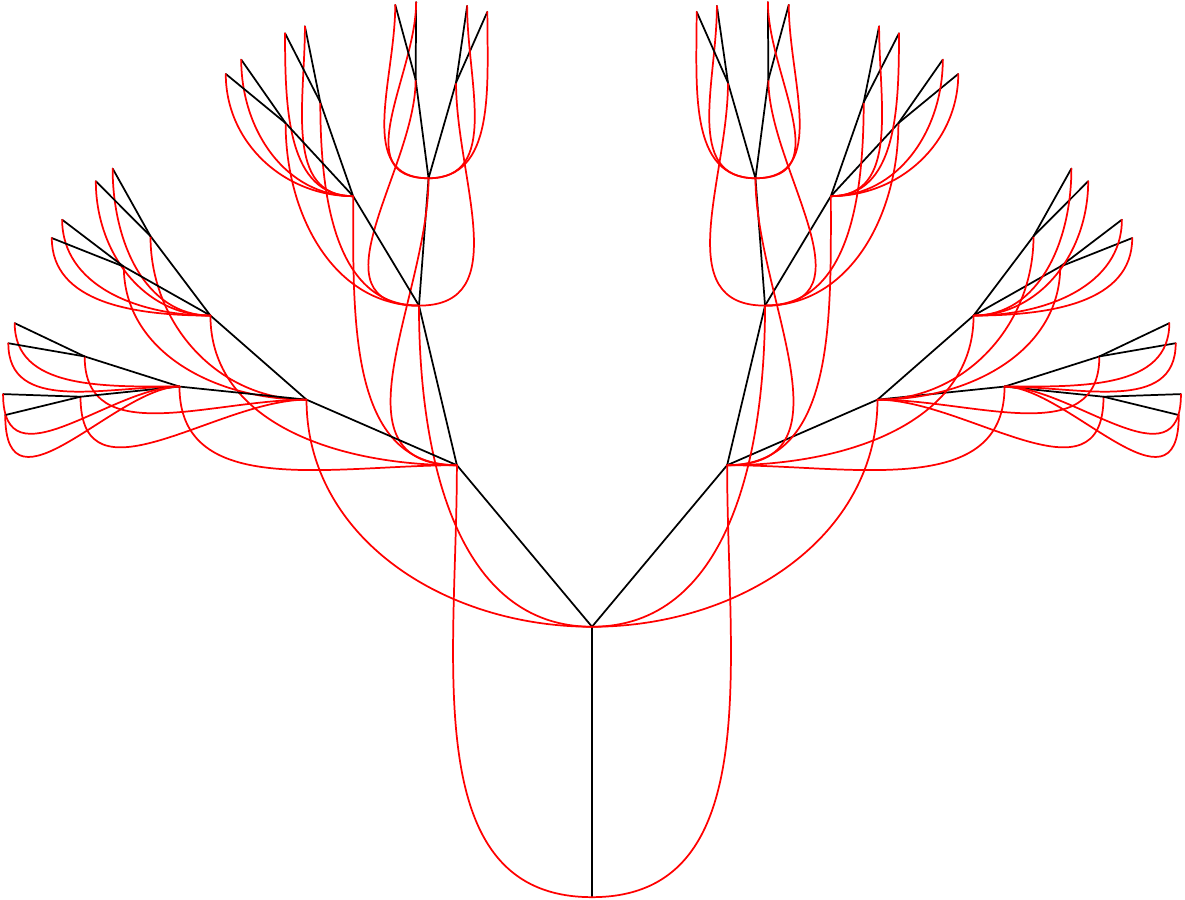}
  \caption{\label{fig:grandfather} The grandfather graph $G_3$. Edges
    of $\bT_3$ are straight black lines. Edges to grandfathers are red
    curves. The distinguished end is the ``down'' direction.}
\end{figure}

The grandfather graph of order $n \geq 3$, $G_n$, is the following
graph (see Figure~\ref{fig:grandfather}). Let $\bT_n$ be the regular
tree of degree $n$. The ends of $\bT_n$ can be identified with the set
of infinite simple paths starting at $o$, an arbitrary distinguished
vertex. Choose a distinguished end. Then each vertex has a unique edge
in the direction of this end. Call the vertex on the other side of
that edge the ``father''. Then each vertex has a unique father, and,
as one can imagine, each vertex has a unique ``grandfather''. The set
of vertices of $G_n$ is identical to that of $\bT_n$. The set of edges
includes the set of edges of $\bT_n$, and in addition an edge between
each vertex and its grandfather.

\subsection{Main result}
\begin{maintheorem}
  \label{thm:main}
  For $n \geq 3$, the grandfather graph $G_n$ is isolated.
\end{maintheorem}
In fact, we show below that the ball of radius one determines the
structure of $G_n$.

We state here without proof that this result can be further extended
to some classes of graphs that are similar to $G_n$. For example, the
product of $G_n$ with any finite graph will also be isolated, as will
``great\textsuperscript{$k$}-grandfather'' graphs.

\section{Proof}
A {\em directed edge} in an undirected graph $G = (V,E)$ is an ordered
pair $(u,w)$ of vertices in $G$ such that $(u,w) \in E$.

Let $(u,w)$ and $(u',w')$ be two directed edges in a graph $G$. We say
that $(u,w)$ and $(u',w')$ are isomorphic if there exists a graph
isomorphism of $G$ that maps $u$ to $u'$ and $w$ to $w'$ (compare to
the notion of ``doubly rooted graphs'' - see,
e.g.,~\cite{aldous2007processes, lyons2013probability}). While all
vertices in a transitive graph are isomorphic, not all directed edges
are necessarily isomorphic.

\begin{figure}[h]
  \centering
  \includegraphics{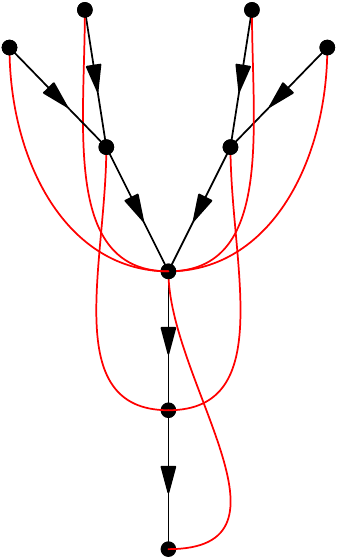}
  \caption{\label{fig:ball} The ball of radius one in the
    grandfather graph $G_3$. The directions and labels of the edges
    can be inferred from the undirected graph.}
\end{figure}

In the grandfather graph $G_n$, $(u,w)$ and $(u',w')$ are isomorphic
if and only if both pairs can be described by the same (ordered)
familial relation: that is, if $w$ is $u$'s father (respectively son /
grandfather / grandson) and $w'$ is $u'$'s father (respectively son /
grandfather / grandson). This is a well-known property of this graph
(see, e.g.,~\cite{benjamini2012ergodic, lyons2013probability}), and is
in fact easy to already see by examining the ball of radius one (see
Figure~\ref{fig:ball}). In this subgraph, the father and the sons can
be distinguished from the grandfather and the grandsons, since their
degrees differ. Furthermore, the father can be distinguished from the
sons, since the father is connected to the sons (he is their
grandfather), while they are not connected to each other. Thus, if $w$
is $u$'s father but $w'$ is not $u'$'s father, there is no graph
isomorphism of $G_n$ that maps $(u,w)$ to $(u',w')$.

We can therefore label each directed edge as a father / son /
grandfather / grandson edge (that is, $(u,w)$ will be a father edge if
$w$ is $u$'s father), and this labeling will be invariant to any
isomorphism of the graph.

This labeling gives rise to an equivalent definition of the
grandfather graph: define a father relation on the $n$-regular tree
$\bT_n$; this is any relation in which each node has a unique father
which is its neighbor in the graph. Then, connect each node to its
grandfather. The choice of a father relation is equivalent to a choice
of end, and hence this also results in the grandfather graph.

Let $H$ be any graph in $\cG$ such that $B_1(H) = B_1(G_n)$. Since
examining the ball of radius one around each vertex is sufficient to
determine this labeling, we can also label the directed edges of $H$
in the same manner, and this labeling will also be invariant to the
isomorphism group of $H$.  We will use this to show that $H$ is
isomorphic to $G_n$, which will prove Theorem~\ref{thm:main}.  

A simple cycle in a graph is a sequence of directed edges $(u_0,w_0),
\ldots, (u_{k-1},w_{k-1})$ such that $w_i = u_{i+1\mod k}$, and each
edge is visited at most once.
\begin{claim}
  \label{clm:no-cycles}
  There are no simple cycles in $H$ which are comprised only of father
  edges and of son edges.
\end{claim}
\begin{proof}
  Assume by contradiction that $(u_0,w_0), \ldots, (u_{k-1},w_{k-1})$
  is a simple cycle comprised only of father edges and and son
  edges. Then all edges are of the same type (i.e., all father edges
  or all son edges): otherwise, there must be in the cycle a father
  edge followed by a son edge, which would make the cycle non-simple,
  since fathers are unique.

  By changing the direction of the cycle we can therefore assume
  without loss of generality that all edges are father edges. Now,
  note that since $B_1(H) = B_1(G_n)$, it also follows that every node
  in $H$ has a unique father and exactly $n-1$ sons. Hence each node
  on the cycle is its own $k$\textsuperscript{th}-order father, and
  each node has $n-2>0$ sons which are not on the cycle. Since the
  father relation is invariant to graph isomorphisms, so is the
  $k$\textsuperscript{th}-order father relation.

  Let $u$ be a vertex on the cycle, and let $v$ be vertex which is not
  on the cycle and is a son of $u$. Then there is no graph isomorphism
  of $H$ that maps $v$ to $u$, since $v$ - unlike $u$ - is not its own
  $k$\textsuperscript{th}-order father. Hence $H$ is not transitive,
  and we have reached a contradiction.
\end{proof}

{\bf Remark. } This claim can also be proved by showing that $H$ is
not unimodular and analyzing the Haar measure of the stabilizers of
the nodes lying on the cycle (see~\cite{lyons2013probability}).
\\
\\
It follows immediately from Claim~\ref{clm:no-cycles} that the
restriction of $H$ to father-son edges is isomorphic to $\bT_n$, the
$n$-regular tree. This restriction is still a connected graph, since
grandfather edges only connect nodes already connected by length two
paths of father edges.

Since $B_1(H) = B_1(G_n)$, the grandfather edges in $H$ are determined
by the father-son relation, and in the same way that they are
determined in $G_n$. Hence $H$ can be constructed by adding
grandfather edges to $\bT_n$, equipped with a father relation. It
follows that $H$ is isomorphic to $G_n$, thus proving Theorem~\ref{thm:main}.

\bibliography{grandfather}

\begin{bibdiv}
\begin{biblist}

\bib{aldous2007processes}{article}{
      author={Aldous, David},
      author={Lyons, Russell},
       title={Processes on unimodular random networks},
        date={2007},
        ISSN={1083-6489},
     journal={Electronic Journal of Probability},
      volume={12},
       pages={no. 54, 1454\ndash 1508},
         url={http://dx.doi.org/10.1214/EJP.v12-463},
      review={\MR{2354165 (2008m:60012)}},
}

\bib{aldous2003objective}{article}{
      author={Aldous, David},
      author={Steele, J.~Michael},
       title={The objective method: Probabilistic combinatorial optimization
  and local weak convergence},
        date={2003},
        ISSN={0938-0396},
     journal={Probability on Discrete Structures (Volume 110 of Encyclopaedia
  of Mathematical Sciences), ed. H. Kesten},
      volume={110},
       pages={1\ndash 72},
}

\bib{benjamini2012ergodic}{article}{
      author={Benjamini, Itai},
      author={Curien, Nicolas},
       title={Ergodic theory on stationary random graphs},
        date={2012},
     journal={Electron. J. Probab},
      volume={17},
      number={93},
       pages={1\ndash 20},
}

\bib{benjamini2001recurrence}{article}{
      author={Benjamini, Itai},
      author={Schramm, Oded},
       title={Recurrence of distributional limits of finite planar graphs},
        date={2001},
     journal={Electronic Journal of Probability},
       pages={1\ndash 13},
}

\bib{champetier2000espace}{article}{
      author={Champetier, Christophe},
       title={L'espace des groupes de type fini},
        date={2000},
     journal={Topology},
      volume={39},
      number={4},
       pages={657\ndash 680},
}

\bib{cornulier2009cantor}{article}{
      author={Cornulier, Yves},
       title={On the cantor-bendixson rank of metabelian groups},
        date={2009},
     journal={arXiv preprint arXiv:0904.4230},
}

\bib{de2007isolated}{article}{
      author={Cornulier, Yves},
      author={Guyot, Luc},
      author={Pitsch, Wolfgang},
       title={On the isolated points in the space of groups},
        date={2007},
     journal={Journal of Algebra},
      volume={307},
      number={1},
       pages={254\ndash 277},
}

\bib{lyons2013probability}{book}{
      author={Lyons, Russell},
      author={Peres, Yuval},
       title={Probability on tress and networks},
        date={2013},
}

\bib{trofimov1985automorphism}{article}{
      author={Trofimov, Vladimir},
       title={Automorphism groups of graphs as topological groups},
        date={1985},
     journal={Matematicheskie Zametki},
      volume={38},
      number={3},
       pages={378\ndash 385},
}

\end{biblist}
\end{bibdiv}
\end{document}